\documentclass[12pt]{amsart}
\usepackage{amsthm,amsmath,amssymb}
\usepackage[utf8]{inputenc}
\DeclareMathOperator{\Aff}{Aff}

\newtheorem{theorem}{Theorem}
\newtheorem{conjecture}{Conjecture}
\theoremstyle{definition}

\begin{document}
\title{Enumeration of strong dichotomy patterns}
\author{Octavio A. Agust\'in-Aquino}
\address
{Instituto de Física y Matemáticas, Universidad Tecnológica de la Mixteca, Carretera a Acatlima km 2.5, Huajuapan de León, Oaxaca, México, C.P. 69000}
\email{octavioalberto@mixteco.utm.mx}

\date{September 23, 2018}

\subjclass[2010]{00A65, 05E18}
\keywords{strong dichotomy pattern, P\'{o}lya-Redfield theory, cyclic sieving}

\maketitle

\begin{abstract}
We apply the version of P\'{o}lya-Redfield theory obtained by White to count patterns with
a given automorphism group to the enumeration of strong dichotomy patterns, that is,
we count bicolor 
 patterns of $\mathbb{Z}_{2k}$ with respect to the action of $\Aff(\mathbb{Z}_{2k})$ and with
trivial isotropy group. As a byproduct, a conjectural instance of phenomenon similar to
cyclic sieving for special cases of these combinatorial objects is proposed.
\end{abstract}

\section{Introduction}

In a short and beautiful paper \cite{dW75}, White proved an analogue of Cauchy-Frobenius-Burnside
lemma tailored for the purpose of counting patterns with a fixed group of automorphisms. Before stating it, 
we warn the reader that we will use the Iverson bracket\footnote{White uses the similar notation $\chi(P)$ in his articles, which was introduced first by Adriano Garsia
in a paper from 1979 according to Knuth \cite{dK92}, although White's paper predates it by four years.} as defined
by Graham, Knuth and Patashnik \cite[p. 24]{GKP94}: if $P$ is a property, then
\[
[P] = \begin{cases}
1, & P\text{ is true},\\
0, & \text{otherwise}.
\end{cases}
\]

\begin{theorem}[D. E. White, 1975]\label{T:Burnsidelike}
Let $S$ be a finite set, $G$ a finite group acting on $S$ and $\Delta$ a system of orbit
representatives for $G$ acting on $S$. Suppose $\{G_{1},\ldots,G_{N}\}$ is a transversal of the orbits
of the subgroups of $G$ under the conjugation action, such that
\[
 |G_{1}|\geq \cdots \geq |G_{N}|.
\]

Given a weight function $w:S\to T$ such that $w(\sigma s) = w(s)$
for all $\sigma \in G$ and all $s\in S$, we have
\[
\sum_{s\in\Delta}w(s)[G_{s}\sim G_{i}] = \sum_{j=1}^{N}b_{i,j}\sum_{s\in S}w(s)[G_{j}s = s],
\]
where $B = (b_{i,j})$ is the inverse of the \emph{table of marks} matrix
\[
M_{i,j} =  \frac{1}{|G_{j}|}\sum_{\sigma\in G}[\sigma G_{i}\sigma^{-1}\subseteq G_{j}].
\]
\end{theorem}

Note that the table of marks matrix is invertible because it is triangular and no element
of the diagonal is $0$.

Let $S= R^{D}$, where both $R$ and $D$ are finite sets. If $G$ acts on $R$ (the set of \emph{colors}),
it is well known that the action can be extended to $S$ defining $\sigma\cdot f = f\circ \sigma^{-1}$, where
$\sigma$ is reinterpreted as a member of the permutation group of $D$.
We know that the action of $G_{i}$ on $D$ defines
a set of disjoint orbits
\[
 O_{G_{i}:D}:=\{G_{i}x_{1},\ldots,G_{i}x_{\ell}\}
\]
which is a partition of $D$, so we can define
\[
q_{G_{i}}(d) = \sum_{i=1}^{\ell} [|Gx_{i}|=d].
\]

This allows us to define the \emph{orbit index monomial}
as
\[
P_{i}(z_{1},z_{2},\ldots,z_{|D|}) := \prod_{d\in D} z_{d}^{q_{G_{i}}(d)},
\]
which can be used in a straightforward manner to obtain a pattern inventory polynomial.

\begin{theorem}[D. E. White, 1975]
The pattern inventory polynomial for patterns fixed by the subgroup $G_{i}$ is
\[
Q_{i} = \sum_{j=1}^{N} b_{i,j}P_{j}(z_{1},z_{2},\ldots,z_{|D|}),
\]
where the substitution $y_{i} = \sum_{r\in R} x_{r}^{i}$ is made.
\end{theorem}

We will use White's results (and one further generalization obtained by him that we will discuss later) to count bicolor patterns
under a group $G$ which are of particular interest for
mathematical musicology. The usual choices for $G$ are cyclic, dihedral and
general affine groups, since they model common musically meaningful transformations,
such as transpositions, inversions, retrogradations and others related to
twelve-tone techniques, to name a few (see \cite[Chapter 8]{gM02} for more examples).

One reason to study this kind of combinatorial objects is that they represent rhythmic patterns if they are interpreted as onsets
in a measure (see \cite{hP99} and \cite{HK06} and the references therein for more information).
Another reason is that they can be seen as abstractions of the concepts of consonance and dissonance in Renaissance counterpoint. In particular,
\emph{self-complementary} (that is, those whose complement belongs to its orbit) and \emph{rigid} (which
means that they are invariant only under the identity) patterns, hereafter called \emph{strong}, are known
to be used in both Western and Eastern music \cite[Part VII]{gM02}, and that their combinatorial structure
lead to significant musicological results \cite[Chapter 31]{gM02}. Note in passing that self-complementarity forces the patterns to be
subsets of cardinality $k$ of sets of even cardinality $2k$. In general, \emph{dichotomy} patterns are those of cardinality
$k$ within a set of cardinality $2k$.

In the following section we provide simple examples of White's theory in action to explain the algorithms we use in the
main computations. The results of these calculations appear in Section \ref{S:Calc}, and in the final section we provide some further
comments regarding them.

\section{Two easy examples}

Suppose we color black or white the vertices of a rectangle that is not a square. The group of symmetries acting on
the colorings of the vertices is the Klein four-group
\[
V = \langle a,b|a^{2}=b^{2}=(ab)^{2}=e\rangle.
\]

We will find the patterns that are invariant under $G_{1}=V$, $G_{2}=\langle a\rangle$, $G_{3}=\langle b\rangle$,
$G_{4} =\langle ab\rangle$ and $G_{5}=\langle e\rangle$ using White's formulas. Since all the proper subgroups of $V$ are normal,
we easily calculate the table of marks matrix
\[
M_{V} = \begin{pmatrix}
1 & 0 & 0 & 0 & 0\\
1 & 2 & 0 & 0 & 0\\
1 & 0 & 2 & 0 & 0\\
1 & 0 & 0 & 2 & 0\\
1 & 2 & 2 & 2 & 4
\end{pmatrix},
\]
whose inverse is
\[
B = \begin{pmatrix}
1 & 0 & 0 & 0 & 0\\
-\tfrac{1}{2} & \tfrac{1}{2} & 0 & 0 & 0\\
-\tfrac{1}{2} & 0 & \tfrac{1}{2} & 0 & 0\\
-\tfrac{1}{2} & 0 & 0 & \tfrac{1}{2} & 0\\
\tfrac{1}{2} & -\tfrac{1}{4} & -\tfrac{1}{4} & -\tfrac{1}{4} & \tfrac{1}{4}
\end{pmatrix}.
\]

It is illustrative to make explicit the formula of Theorem \ref{T:Burnsidelike} for this simple example. Let us code the
colorings of the vertices with the strings $u_{1}u_{2}u_{3}u_{4}$ over the alphabet $\{n,b\}$, using the clockwise order and
beginning from the upper left corner. Then
\[
\Delta = \{nnnn,nnnb,nnbb,nbbn,nbnb,nbbb,bbbb\}.
\]

For $G_{1}$ the formula trivially asserts that the only patterns invariant under the full group are
the monochromatic ones. For $G_{2}$ we have
\begin{multline*}
-\tfrac{1}{2}(w(nnnn)+w(bbbb))\\
+\tfrac{1}{2}(w(nnnn)+w(nnbb)+w(bbnn)+w(bbbb))\\
= \tfrac{1}{2}(w(bbnn)+w(nnbb)) = w(nnbb)
\end{multline*}
because $w(bbnn)=w(nnbb)$, by hypothesis. The colorings $bbnn$ and $nnbb$ are precisely those who represent the only pattern which is invariant
under the reflection with vertical axis. The cases of $G_{3}$ and $G_{4}$ are analogous. Finally, the case of
the trivial subgroup is more interesting:
\begin{multline*}
\tfrac{1}{2}(w(nnnn)+w(bbbb))\\
-\tfrac{1}{4}(w(nnnn)+w(nnbb)+w(bbnn)+w(bbbb))\\
-\tfrac{1}{4}(w(nnnn)+w(nbbn)+w(bnnb)+w(bbbb))\\
-\tfrac{1}{4}(w(nnnn)+w(nbnb)+w(bnbn)+w(bbbb))+\tfrac{1}{4}\sum_{\text{all strings}} w(s) =\\
\tfrac{1}{4}(w(nnnb)+w(nnbn)+w(nbnn)+w(nbbb)\\+w(bnnn)+w(bnbb)+w(bbnb)+w(bbbn))\\
=w(nnnb)+w(nbbb), 	
\end{multline*}
and it informs us of the two patterns that are invariant under the action of the trivial subgroup only;
they are precisely those with only one black or only one white vertex.

Let us confirm the former using the orbit index polynomials for each subgroup. For $G_{1}$, we have only one
orbit of four elements, thus
\[
P_{1} = z_{4}.
\]

The orbits defined by $G_{2}$, $G_{3}$ and $G_{4}$ are all of cardinality two, thus
\[
P_{2}=P_{3}=P_{4} = z_{2}^{2}.
\]

Finally, there are four orbits of cardinality one for the trivial subgroup, hence
\[
P_{5} = z_{1}^{4}.
\]

Using these polynomials, we can calculate all the pattern inventories at once:
\begin{align*}
\begin{pmatrix}
Q_{1}\\
Q_{2}\\
Q_{3}\\
Q_{4}\\
Q_{5}
\end{pmatrix}
&=
\begin{pmatrix}
1 & 0 & 0 & 0 & 0\\
-\tfrac{1}{2} & \tfrac{1}{2} & 0 & 0 & 0\\
-\tfrac{1}{2} & 0 & \tfrac{1}{2} & 0 & 0\\
-\tfrac{1}{2} & 0 & 0 & \tfrac{1}{2} & 0\\
\tfrac{1}{2} & -\tfrac{1}{4} & -\tfrac{1}{4} & -\tfrac{1}{4} & \tfrac{1}{4}
\end{pmatrix}
\begin{pmatrix}
z_{4} \\ z_{2}^{2} \\ z_{2}^{2} \\ z_{2}^{2} \\ z_{1}^{4}
\end{pmatrix}
\\
&=\begin{pmatrix}
z_{4}\\
\tfrac{1}{2}z_{2}^{2}-\tfrac{1}{2}z_{4}\\
\tfrac{1}{2}z_{2}^{2}-\tfrac{1}{2}z_{4}\\
\tfrac{1}{2}z_{2}^{2}-\tfrac{1}{2}z_{4}\\
\tfrac{1}{2}z_{4}-\tfrac{3}{4}z_{2}^{2}+\tfrac{1}{4}z_{1}^{4}
\end{pmatrix}.
\end{align*}

Upon the substitution $z_{i} = 1+x^{i}$, that allows us to count
the number of bicolor patterns according to the number of black elements (say), we find
\[
\begin{pmatrix}
Q_{1}\\
Q_{2}\\
Q_{3}\\
Q_{4}\\
Q_{5}
\end{pmatrix}
=
\begin{pmatrix}
1+x^{4}\\
x^{2}\\
x^{2}\\
x^{2}\\
x+x^{3}
\end{pmatrix}.
\]

We proceed now with a more complicated example that introduces the group we will use in our main
computation. Define
\[
\Aff(\mathbb{Z}_{2k}) = \mathbb{Z}/2k\mathbb{Z}\ltimes \mathbb{Z}/2k\mathbb{Z}^{\times}.
\]

Denote an element $(u,v)\in\Aff(\mathbb{Z}_{2k})$ by $e^{u}.v$. The action of
$\Aff(\mathbb{Z}_{2k})$ on $\mathbb{Z}/2k\mathbb{Z}$ is given by
\[
 e^{u}.v(x) = vx+u.
\]

Let us compute the number of patterns of the action of $\Aff(\mathbb{Z}_{6})$. We have
the following sequence of normal subgroups,
\begin{multline*}
G_{1} = \langle e^{1}.1, e^{0}.5\rangle,
G_{2} = \langle e^{4}.1, e^{5}.5\rangle,\\
G_{3} = \langle e^{1}.1, e^{0}.5\rangle,
G_{4} = \langle e^{1}.1\rangle,\\
G_{5} = \langle e^{3}.1, e^{0}.5\rangle,
G_{6} = \langle e^{2}.1\rangle,\\
G_{7} = \langle e^{5}.5\rangle,
G_{8} = \langle e^{0}.5\rangle,\\
G_{9} = \langle e^{3}.1\rangle,
G_{10} = \{e^{0}.1\},
\end{multline*}

The computation of the table of marks matrix is not as direct as before, in part because the
subgroups $G_{5}$, $G_{7}$ and $G_{9}$ are not normal. But using GAP \cite{GAP4} we readily find
\[
M = \begin{pmatrix}
1 &0 &0 &0 &0 &0 &0 &0 &0 &0\\
1 &2 &0 &0 &0 &0 &0 &0 &0 &0\\
1 &0 &2 &0 &0 &0 &0 &0 &0 &0\\
1 &0 &0 &2 &0 &0 &0 &0 &0 &0\\
1 &0 &0 &0 &1 &0 &0 &0 &0 &0\\
1 &2 &2 &2 &0 &4 &0 &0 &0 &0\\
1 &2 &0 &0 &1 &0 &2 &0 &0 &0\\
1 &0 &2 &0 &1 &0 &0 &2 &0 &0\\
1 &0 &0 &2 &3 &0 &0 &0 &6 &0\\
1 &2 &2 &2 &3 &4 &6 &6 &6 &12
\end{pmatrix}
\]
whose inverse is
\[
B = \begin{pmatrix}
\tfrac{1}{12} &0 &0 &0 &0 &0 &0 &0 &0 &0\\
-\tfrac{1}{12} &\frac{1}{6} &0 &0 &0 &0 &0 &0 &0 &0\\
-\tfrac{1}{4} &0 &\tfrac{1}{2} &0 &0 &0 &0 &0 &0 &0\\
-\tfrac{1}{4} &0 &0 &\tfrac{1}{2} &0 &0 &0 &0 &0 &0\\
-\tfrac{1}{12} &0 &0 &0 &\tfrac{1}{4} &0 &0 &0 &0 &0\\
\tfrac{1}{2} &-\tfrac{1}{2} &-\tfrac{1}{2} &-\tfrac{1}{2} &0 &1 &0 &0 &0 &0\\
\tfrac{1}{12} &-\tfrac{1}{6} &0 &0 &-\tfrac{1}{4} &0 &\tfrac{1}{2} &0 &0 &0\\
\tfrac{1}{4} &0 &-\tfrac{1}{2} &0 &-\tfrac{1}{4} &0 &0 &\tfrac{1}{2} &0 &0\\
\tfrac{1}{4} &0 &0 &-\tfrac{1}{2} &-\tfrac{1}{4} &0 &0 &0 &\tfrac{1}{2} &0\\
-\tfrac{1}{2} &\tfrac{1}{2} &\tfrac{1}{2} &\frac{1}{2} &\tfrac{1}{2} &-1 &-\tfrac{1}{2} &-\tfrac{1}{2} &-\tfrac{1}{2} &1
\end{pmatrix}
\]

The orbit index polynomials are
\begin{multline*}
P_{1} = z_{6}, P_{2} = z_{6},
P_{3} = z_{3}^{2}, P_{4} = z_{6},
P_{5} = z_{2}z_{4}, \\
P_{6} = z_{3}^{2},
P_{7} = z_{2}^{3}, P_{8} = z_{1}^{2}z_{2}^{2},
P_{9} = z_{2}^{3}, P_{10} = z_{1}^{6}
\end{multline*}
whence
\[
Q(z_{1},\ldots,z_{6})=BP = \begin{pmatrix}
z_{6}\\
0\\
\tfrac{1}{2}z_{3}^{2}-\tfrac{1}{2}z_{6}\\
0\\
z_{2}z_{4}-z_{6}\\
0\\
-\tfrac{1}{2}z_{2}^{3}-\tfrac{1}{2}z_{2}z_{4}\\
\tfrac{1}{2}z_{6}-\tfrac{1}{2}z_{2}z_{4}-\tfrac{1}{2}z_{3}^{2}+\tfrac{1}{2}z_{1}^{2}z_{2}^{2}\\
\tfrac{1}{3}z_{6}-\tfrac{1}{2}z_{2}z_{4}+\tfrac{1}{6}z_{2}^{3}\\
-\tfrac{1}{6}z_{6}+\tfrac{1}{2}z_{2}z_{4}+\tfrac{1}{6}z_{3}^{2}-\tfrac{1}{3}z_{2}^{3}-\tfrac{1}{4}z_{1}^{2}z_{2}^{2}+\tfrac{1}{12}z_{1}^{6}
\end{pmatrix}
\]
thus
\[
Q(1+x,\ldots,1+x^{6}) = \begin{pmatrix}
x^{6}+1\\
0\\
x^{3}\\
0\\
x^{4}+x^{2}\\
0\\
x^{4}+x^{2}\\
x^{5}+x^{4}+x^{3}+x^{2}+x\\
0\\
x^{3}
\end{pmatrix}.
\]

It is interesting to learn that there are no patterns that are exclusively invariant under the subgroups
generated, respectively, by the translations $e^{1}.1$, $e^{2}.1$, $e^{3}.1$. In other words:
arithmetic progressions with common difference $1$, $2$ and $3$ are invariant under symmetries
that are not translations.

\section{Main calculations}\label{S:Calc}

Denoting by $\mathcal{D}$ the set of dichotomies, and by $\mathcal{S}$ and $\mathcal{R}$ the subsets of the self-complementary and rigid dichotomies (respectively), we know by the principle of inclusion and exclusion (PIE) that
\[
|\mathcal{D}| \geq |\mathcal{S}\cup \mathcal{R}| = |\mathcal{S}|+|\mathcal{R}|-|\mathcal{S}\cap \mathcal{R}|
\]
where $|\mathcal{S}\cap \mathcal{R}|$ is precisely the number of strong dichotomies. Hence
\[
 |\mathcal{S}\cap \mathcal{R}| \geq |\mathcal{S}|+|\mathcal{R}|-|\mathcal{D}|.
\]

We can calculate $|\mathcal{D}|$ and $|\mathcal{S}|$ with the classical P\'olya-Redfield theory, and
$|\mathcal{R}|$ with White's formulas, so we may expect this inequality to provide reasonably good bounds on the number
of strong dichotomies. But, unfortunately, in general it does not, as we can readily see in Table \ref{T:PRW}, since
many of them are negative.

However, not everything is lost. After examining the cases when the PIE yields a nontrivial bound, we discover that
this happens when $k$ is a power of a prime and, more importantly, the value of $|Q_{1}(-1)|$ coincides with
the number of strong dichotomy patterns calculated by direct construction for these cases (see \cite{oA11}).
On the other hand, it is known that the
classical pattern inventory polynomials of the P\'olya-Redfield theory exhibit a form of the \emph{cyclic sieving phenomenon} \cite[Corollary 6.2]{RSW04},
which means that if $p(x)$ is the generating function of the number
of patterns according to its number of black elements, then $p(-1)$ yields the number of self-complementary patterns.

Since the polynomials for White's formulas do not count cycles but orbits, in general they fail
to cyclically sieve patterns, but we may expect it to work when $\mathbb{Z}_{2k}^{\times}$ is cyclic. Indeed,
if the group of units is generated by a single element, it is plausible to think that all the orbits are cycles of $e^{1}.1$
and a generator of $\mathbb{Z}_{2k}^{\times}$. Furthermore, it is a well-known fact that the group
of units of $\mathbb{Z}_{n}$ is cyclic precisely
when $n=1,2,p^{k}$, where $p$ is a prime number. This discussion, however, is not a full proof, so we
formalize it as a conjecture. Hopefully, it will be proved soon.

\begin{conjecture}\label{C:CasoFacil}
Let $G_{N}=\Aff(\mathbb{Z}_{2k})$ and $\{G_{i}\}$ be a set of representatives of the orbits of the
conjugation action such that $|G_{N}|\geq \cdots \geq |G_{1}|$ and let $B=(b_{i,j})$ be the inverse of
its table of marks. If $k$ is equal to $1$, $2$ or a power
of an odd prime number, then the pattern inventory polynomial for bicolor patterns fixed by the subgroup $G_{i}$ 
\[
Q_{i} = \sum_{j=1}^{N} b_{i,j}P_{j}(1+x,1+x^{2},\ldots,1+x^{|D|}),
\]
is such that $|Q_{i}(-1)|$ counts the number of self-complementary dichoto\-mies with automorphism
group $G_{i}$. In particular, $|Q_{1}(-1)|$ counts the number of strong dichotomies.
\end{conjecture}

\begin{table}
\begin{center}
\begin{tabular}{|c|c|c|c|c|c|}
\hline
$2k$ & $|\mathcal{D}|$  & $|\mathcal{S}|$ & $|\mathcal{R}|$ & PIE bound  & $|Q_{1}(-1)|$\\
\hline
$2$ & $1$ & $1$ & $1$ & $1$ & $1$\\
$4$ & $2$ & $2$ & $0$ & $0$ & $0$\\
$6$ & $3$ & $3$ & $0$ & $0$ & $1$\\
$8$ & $6$ & $4$ & $1$ & $-1$ & $1$\\
$10$ & $9$ & $7$ & $5$ & $3$ & $3$\\
$12$ & $34$ & $18$ & $10$ & $-6$ & $4$\\
$14$ & $47$ & $15$ & $37$ & $5$ & $9$\\
$16$ & $129$ & $21$ & $83$ & $-25$ & $1$\\
$18$ & $471$ & $55$ & $436$ & $20$ & $40$\\
$20$ & $1280$ & $134$ & $1052$ & $-94$ & $66$\\
$22$ & $3235$ & $115$ & $3181$ & $61$ & $105$\\
$24$ & $15008$ & $440$ & $13331$ & $-1237$ & $33$\\
$26$ & $33429$ & $385$ & $33253$ & $209$ & $355$\\
$28$ & $121466$ & $1194$ & $117422$ & $-2850$ & $886$\\
$30$ & $648819$ & $3365$ & $643901$ & $-1153$ & $3007$\\
$32$ & $1182781$ & $2189$ & $1165498$ & $-15094$ & $1432$\\
$34$ & $4290533$ & $4375$ & $4288913$ & $2755$ & $4305$\\
$36$ & $21082620$ & $18404$ & $20933318$ & $-130898$ & $15518$\\
$38$ & $51677171$ & $15347$ & $51671611$ & $9787$ & $15267$\\
$40$ & $215804540$ & $49684$ & $214972319$ & $-782537$ & $25659$\\
$42$ & $1068159497$ & $133285$ & $1067785287$ & $-240925$ & $130839$\\
$44$ & $2392981542$ & $171662$ & $2389064994$ & $-3744886$ & $155346$\\
$46$ & $8135833183$ & $198943$ & $8135769049$ & $134809$ & $198753$\\
$48$ & $42007923187$ & $786707$ & $41970277573$ & $-36858907$ & $643019$\\
$50$ & $126410742103$ & $872893$ & $126410471144$ & $601934$ & $871992$\\
\hline
\end{tabular}
\end{center}
\caption{Summary of the information that can be obtained via the classical P\'olya-Redfield theory and
White's extension, for $1\leq k \leq 25$.}
\label{T:PRW}
\end{table}

For the general case, we can use another formula of White \cite{dW75b}. Now we need
to consider the swapping action on the colors of the patterns simultaneously with that of the affine group,
so we see $G = \Aff(\mathbb{Z}_{2k})\times \mathbb{Z}_{2}$ as acting
both in $R$ and $D$, according to
\[
(\sigma,\tau)\cdot r = \tau\cdot r\quad\text{and}\quad (\sigma,\tau)\cdot d = \sigma \cdot d;
\]
hence $G$ acts doubly on $R^{D}$ in the following manner
\[
 g\cdot f = g\circ f \circ g^{-1}.
\]

We provide a quick sketch of White's reasoning to obtain the counting formula, in part because our problem's
conditions lead to a simpler statement
and in part because his original paper has some minor (but misleading) typographical errors.

In order to
apply Theorem \ref{T:Burnsidelike}, we must characterize
first the subgroups $H\subseteq G_{f}$ for a pattern $f$. If $H$ leaves the pattern $f$ invariant, it means
that $f$ sends an element of $O_{H:D}$ to an element of $O_{H:R}$.

Thus, let $B\in O_{H:D}$ and $f(B)=C\in O_{H:R}$. Taking arbitrary elements $b\in B$
and $c\in C$, we deduce that $f$ must be defined by $f(\gamma_{1} b) = \gamma_{1} c$. This relation, however,
might not be functional, for it may happen that $f(\gamma_{1}b)\neq  f(\gamma_{2}b)$ when $\gamma_{1}b=
\gamma_{2}b$, unless $\gamma_{1}c=\gamma_{2}c$, or $\gamma_{1}^{-1}\gamma_{2}\in H_{c}$. In other words, if
the function is well defined then
\begin{equation}\label{E:Condicion}
 \gamma_{1}^{-1}\gamma_{2}\in H_{b}\quad\text{implies that}\quad \gamma_{1}^{-1}\gamma_{2}\in H_{c}
\end{equation}
or, equivalently,
\[
 H_{b}\subseteq H_{c}
\]
(note that these isotropy groups are relative to $H$). To check that \eqref{E:Condicion} is also
sufficient is direct, like the fact that the election of $b$ is irrelevant.

We have
\[
\sum_{f\in S}w(f)[G_{j}f = f] = \sum_{\hat{f}\in O_{H:R}^{O_{H:D}}}\prod_{B\in O_{H:D}} \sum_{j=0}^{|\hat{f}(B)|-1}[H_{b}\subseteq H_{\tau_{j}c}] w(f)
\]
and, reorganizing the terms (in what White calls \emph{sum-product interchange}), we get
\begin{equation}
\label{E:Circ1}
\sum_{f\in S}w(f)[G_{j}f = f] = \prod_{B\in O_{H:D}}\sum_{C\in O_{H:R}} \sum_{j=0}^{|C|-1} [H_{b}\subseteq H_{\tau_{j}c}] w(f)
\end{equation}
where
\begin{equation}
\label{E:Circ2}
w(f) = \prod_{i=0}^{|B|-1}x_{\tau_{j}c}.
\end{equation}

For bicolor patterns we have $C=\{0,1\}$, therefore the invariance under the action of the whole group
reduces the weights $w(f)$ to the following choices:
\[
w(f) = \begin{cases}
x_{0}^{|B|/2}x_{1}^{|B|/2}, & H\text{ swaps colors},\\
x_{0}^{|B|}=x_{1}^{|B|}, & \text{otherwise}.
\end{cases}
\]

Thus we can reuse the previous algorithm that involves the inverse of the table of marks matrix, but with the larger group
$\Aff(\mathbb{Z}_{2k})\times \mathbb{Z}_{2}$ and
calculating the corresponding vector of polynomials with \eqref{E:Circ1} and \eqref{E:Circ2}. The only
remaining detail is that no longer we may read the total number of strong dichotomy patterns in a single entry of the output vector, for such patterns
have automorphism groups of cardinality two; namely, the identity $(e^{0}.1,0)$ and the color swap $(e^{0}.1,1)$ composed
with a unique symmetry of $\Aff(\mathbb{Z}_{2k})$, which is called the \emph{polarity} of the pattern. Hence, we
gain a feature and not an inconvenience, for now we can know the number of strong dichotomy patterns for each polarity.

The first case that is not covered by Conjecture \ref{C:CasoFacil} is $n=8$, but the table of marks matrix is of size $148
\times 148$, so we will not display it here. Let us simply state that there is only one strong dichotomy, whose polarity is
$e^{5}.-1$. In Table \ref{T:PRWH} we summarize the information that can be calculated with this algorithm up to $2k=48$.

\begin{table}
\begin{center}
\begin{tabular}{|c|c|}
\hline
$2k$ & $|\mathcal{S}\cap \mathcal{R}|$\\
\hline
$6$ & $1$\\
$8$ & $1$\\
$12$ & $2+4=6$\\
$16$ & $1+14=15$\\
$20$ & $3+6+54+27=90$\\
$24$ & $14+54+63+228=359$\\
$28$ & $38+76+326+652=1092$\\
$32$ & $120+2032=2152$ \\
$36$ & $560+1120+5382+10764=17826$ \\
$40$ & $1572+6357+8100+32520=48549$ \\
$42$ & $3936+12135+28320+86448=130839$ \\
$44$ & $4662+9324+52278+104556=170820$ \\
$48$ & $21435+65040+172410+521760=780645$\\
\hline
\end{tabular}
\end{center}
\caption{Summary of the information that can be obtained via White's extension of P\'olya-Redfield theory for strong
dichotomy patterns and selected values of $k$. The totals
of strong dichotomies are displayed as sums, where each summand represents the number of patterns with a specific polarity. Thus,
the number of summands on each row is the number of polarities.}
\label{T:PRWH}
\end{table}

\section{Concluding remarks}

The enumerations of strong dichotomies done here coincide with the explicit ones performed in \cite{oA11} and subsequent
verifications done by the author, with a variation of the original algorithm presented in \cite{oA11}. It is interesting
to note that Conjecture \ref{C:CasoFacil} is of practical interest, since it significantly simplifies the computation
of the table of marks: we should consider that
the volume of calculations is exacerbated when we have to calculate with the product $\Aff(\mathbb{Z}_{2k})\times \mathbb{Z}_{2}$;
its table of marks can be much bigger that the one of its largest factor.

Harald Fripertinger noted in a personal communication with the author that the number of self-complementary patterns $|\mathcal{S}|$
seems to approach asymptotically to
the number of the strong ones (or, equivalently, that the vast majority of dichotomies is rigid). In particular,
$|\mathcal{S}|$ provides a direct and fast way (it does not
require to compute the table of marks) to determine a very good upper bound for the
number of strong patterns, a useful fact in order to partially validate the exact (but lengthy) calculations.


\bibliographystyle{amsplain}
\bibliography{reporte15}

\end{document}